\documentclass[letterpaper, 10 pt, conference]{ieeeconf}  

\IEEEoverridecommandlockouts                              

\usepackage{tikz}
\usepackage{graphics} 
\usepackage{epsfig} 
\usepackage{amsmath} 
\usepackage{amssymb}  
\usepackage{url,algorithm}
\usepackage{multirow}
\usepackage{color}

\let\labelindent\relax
\usepackage{enumitem}
\renewcommand{\theenumi}{\arabic{enumi}}

\renewcommand{\theenumii}{\arabic{enumii}}

\renewcommand{\theenumiii}{\arabic{enumiii}}

\renewcommand{\theenumiv}{\arabic{enumiv}}

{\normalfont}{\normalfont}

\newcommand{\Nu}{N_u}
\newcommand{\Nc}{\Nu}
\newcommand{\Np}{N_p}
\newcommand{\Ts}{T_{\mathrm{s}}}

\newcommand{\nin}{n_u}
\newcommand{\nstate}{n_x}
\newcommand{\nout}{n_y}

\newcommand{\parvec}{\theta}
\newcommand{\setparvec}{\Theta}
\newcommand{\sizeparvec}{n_\parvec}
\newcommand{\pref}{\pi}
\newcommand{\rr}{{\mathbb R}}
\newcommand{\ba}[1]{\begin{array}{#1}}
	\newcommand{\ea}{\end{array}}
\newcommand{\st}{\mathop{\rm s.t.}\nolimits}

\title{\LARGE \bf
Preference-based MPC calibration
}

\author{ Mengjia Zhu, Alberto Bemporad, Dario Piga
	\thanks{M. Zhu and A. Bemporad are with IMT School for Advanced Studies Lucca,  Lucca, Italy. 
		{\tt\small mengjia.zhu@imtlucca.it;  alberto.bemporad@imtlucca.it}}   %
			\thanks{D. Piga is with IDSIA Dalle Molle Institute for Artificial Intelligence,  SUPSI-USI, Lugano, Switzerland. {\tt\small dario.piga@supsi.ch} } 
			\thanks {This paper was partially supported by the Italian Ministry of University and Research under the PRIN'17 project "Data-driven learning of constrained control systems" , contract no. 2017J89ARP.}
}

\begin{document}

\maketitle
\thispagestyle{empty}
\pagestyle{empty}

\begin{abstract}	
Automating the calibration of the parameters of a control policy
by means of global optimization requires quantifying a 
closed-loop performance function. As this can be impractical in many situations, in this paper we suggest a semi-automated calibration approach that requires instead a human calibrator to express a \emph{preference} on whether a certain control policy is ``better'' than another one, therefore eliminating the need of an explicit performance index. In particular, we focus our attention on  semi-automated calibration of \emph{Model Predictive Controllers} (MPCs), for which we
attempt computing the set of best calibration parameters by employing the recently-developed active preference-based optimization algorithm GLISp. Based
on the preferences expressed by the human operator, GLISp learns a surrogate
of the underlying closed-loop performance index that the calibrator (unconsciously) uses and proposes, iteratively, a new set of calibration parameters to him or her for testing and for comparison against previous experimental results. The resulting semi-automated calibration procedure is tested on two case studies, showing the capabilities of the approach in achieving near-optimal performance within a limited  number of experiments.
\end{abstract}

\section{Introduction}
The design of \emph{Model Predictive Controllers}  (MPC) typically requires the tuning of several parameters such as the prediction and control horizons, the weight matrices defining the cost function, various numerical solver tolerances, etc. A calibrator typically adjusts these knobs  based on  experience and trial-and-error, until the closed-loop system behaves as he or she desires. Such a tuning process thus requires skilled calibrators, domain knowledge, and it can be costly and time consuming. 

To automate the tuning process, usually a figure of merit is defined to quantitatively assess closed-loop performance, and experiment-driven optimization algorithms are usually adopted to  find  near-optimal   MPC parameters. In particular, black-box global optimization methods based on surrogate functions,
like \emph{Bayesian Optimization} (BO) \cite{brochu2010} and the recently introduced GLIS (GLobal optimization based on Inverse distance weighting and radial basis function Surrogates) algorithm~\cite{Bem20}, have been recently applied for MPC parameter tuning~\cite{FoPiBe20,LuFoCo20},   choice of the  MPC predictive model~\cite{piga2019performance,bansal2017goal} 
and also in other control engineering problems and applications such as PID and state-feedback control tuning~\cite{DBLP:journals/corr/abs-1906-12086,marco2016automatic},
 position and force control   in robot  manipulators~\cite{roveda2020two,driess2017constrained},  and control of mobile robots and quadrotors~\cite{calandra2014bayesian,berkenkamp2016safe}, just to cite a few. 


The aforementioned optimization approaches iteratively suggest new parameters to be tested based on a surrogate function, which is estimated from closed-loop performance figures gathered from  previous  experiments. An exploration term takes care of sufficiently covering the set of parameters to search, thus avoiding the solver to be trapped in a local minimum.  

However, in order to use such methods for MPC calibration, it is essential to have a well-defined performance index that captures the desired closed-loop behavior
of the system. Unfortunately, in many practical control  applications a performance  is usually formulated based on multiple criteria, and thus it is difficult for a calibrator to formally define and  quantify objectively the scoring function.  On the other hand, it is usually easier for a calibrator to express a preference (such as ``controller A is better than controller B'') between the outcome of two experiments. 

\subsection{Contribution}
Motivated by the above consideration, this paper proposes a novel  calibration approach to tune the MPC parameters based on pairwise preferences between experiment outcomes. The approach is semi-automated in that the control parameters to be tested
are selected automatically by an algorithm while closed-loop performance is assessed manually by the calibrator, that therefore is no longer required to formulate a performance index upfront.

The proposed approach for preference-based  MPC calibration relies on the  derivative-free
global optimization algorithm  recently developed in~\cite{BePi20},  which iteratively proposes a new comparison to the calibrator  to  make,  based  on  actively  learning  a  surrogate  of  the latent  objective function from past sampled decision vectors and pairwise preferences. Preference-based optimization method is becoming popular in the field of reinforcement learning (RL)~\cite{christiano2017deep} and a comprehensive review is presented in the survey paper~\cite{wirth2017survey}. Inverse (reinforcement) learning and semi-supervised machine learning approaches are also closely related to preference-based learning, which has been applied to parameter tuning, such as for autonomous driving ~\cite{menner2019inverse,rosbach2019driving,wang2017driving}. The main concern of the inverse-learning approach is that it can not improve upon the demonstrations by the expert. The algorithm implemented in this paper~\cite{BePi20}, called GLISp (preference-based GLIS) proposed a different acquisition function form that balances the trade-off between exploitation and exploration. GLISp has been shown to be very efficient in terms of number of experiments and comparisons  required to compute the global optimum, which is a major drawback of many preference-based RL methods~\cite{wirth2017survey}. 

We show the efficiency of the proposed preference-based
MPC calibration in two case studies. The first one considers the   control of a \emph{Continuous Stirring Tank Reactor}  (CSTR), while the second one is related to autonomous driving of a vehicle with obstacle avoidance. In both cases, overall satisfactory  performance is achieved  within a relatively small number of
closed-loop experiments, without the hard and time-consuming need to specify a quantitative scoring function driving a fully-automated MPC calibration.   


The rest of this paper is organized as follows. Section~\ref{sec:problem} describes the MPC calibration problem. The preference-based tuning approach based on GLISp  is presented  in Section~\ref{sec:pre_tunning}. The application of the proposed approach for two case studies are  discussed  in Section~\ref{sec:case_study}. Finally, conclusions and directions for future research  are drawn in Section \ref{sec:conclusion}.


\section{Problem description}\label{sec:problem}
Let us consider the problem of controlling a nonlinear multi-input multi-output  system described by the continuous-time state-space representation:
\begin{equation} 
	\begin{split}
	\dot x &= f(x,u) \\
	y &= g(x,u),
	\end{split}%
	\label{mathcalS}%
\end{equation}
where  $x\in \mathbb{R}^{\nstate}$ and $\dot x\in \mathbb{R}^{\nstate}$ are
the state vector and its time derivative, respectively;   $u\in \mathbb{R}^{\nin}$ is the  control input; $y\in \mathbb{R}^{\nout}$ is the vector  of controlled outputs; and $f: \mathbb{R}^{\nstate} \times \mathbb{R}^{\nin} \to\mathbb{R}^{\nstate}$ and $g: \mathbb{R}^{\nstate} \times \mathbb{R}^{\nin} \to\mathbb{R}^{\nout}$ are the state and output mappings, respectively. 

A popular strategy to achieve a reference-tracking objective for a system described by~\eqref{mathcalS} under input and output constraints, is linear-time varying MPC, a.k.a. real-time iteration scheme~\cite{GZQBD16,diehl2005realtimeIter}. The MPC is designed based on the following predictive model obtained via linearization  of~\eqref{mathcalS} around  a nominal trajectory $\bar{x}_k$, $\bar{u}_k$, $\bar{y}_k$ and discretization with sampling time $\Ts$, resulting
in the prediction model
\begin{equation}  
	\begin{split}
	 \tilde{x}_{k+1}&=A_k\tilde{x}_k+B_k \tilde{u}_k  \\ 
	 \tilde{y}_k&=C_k \tilde{x}_k+D_k \tilde{u}_k,  
	\label{eq:MPCsys}
	\end{split}
\end{equation}
where $\tilde{x}_k = x_k - \bar{x}_k$, $\tilde{u}_k = u_k - \bar{u}_k$ and $\tilde{y}_k = y_k - \bar{y}_k$.

At each sampling time $t$, the MPC action $u_{t|t}$ to apply to the system is computed by solving the Quadratic Programming (QP) problem:
\begin{equation} \label{eq:MPC}
	\begin{split}
	&  \min_{\left\{u_{t+k|t}\right\}_{k=0}^{\Nu-1}, \varepsilon}  \sum_{k=0}^{\Np-1} \left\|y_{t+k|t}-y^{\rm{ref}}_{t+k}\right\|_ {Q_y}^2 \\ 
	& \qquad \qquad \qquad+\sum_{k=0}^{\Np-1} \left\|u_{t+k|t}-u^{\rm ref}_{t+k}\right\|_{Q_u}^2  \!\!\\
	&\qquad  \qquad \qquad + \sum_{k=0}^{\Np-1} \left\| \Delta u_{t+k|t} \right\|_{Q_{\Delta u}}^2  + Q_\varepsilon \left\| \varepsilon\right\|^2 \\
	\end{split}
\end{equation}
  $\mathrm{s.t. \ }$ model equation~\eqref{eq:MPCsys} and the following constraints
\begin{equation}  \label{eq:MPC_constraints}
	\begin{split}
	& \quad y_{\mathrm{min}} \!-\! \varepsilon\leq y_{t+k|t} \leq  y_{\mathrm{max}}+\varepsilon, \  k=1,\ldots,\Np  \\
	& \quad   u_{\mathrm{min}}  \leq u_{t+k|t} \leq  u_{\mathrm{max}}, \  k=1,\ldots,\Np \\
	& \quad  \Delta u_{\mathrm{min}}   \leq \Delta u_{t+k|t}, \ \    k=1,\ldots,\Np   \\
	& \quad    \Delta u_{t+k|t} \leq  \Delta u_{\mathrm{max}}, \ \  k=1,\ldots,\Np \\ 
	& \quad  u_{t+\Nu+j|t}= u_{t+\Nu|t},  \ \  j=1,\ldots, \Np-\Nu, 
	\end{split}
\end{equation}
where $\left\|\nu\right\|^2_Q$ is the weighted squared norm, \emph{i.e.}, $\nu^{\top}Q\nu$; $u_{\rm ref}$ and $y_{\rm ref}$ are the input and output references, respectively;  and  $\Delta u_{t+k|t} =  u_{t+k|t}- u_{t+k-1|t}$. 

Several tuning parameters appear in the MPC problem~\eqref{eq:MPC} -~\eqref{eq:MPC_constraints} and must be tuned, such as:
\begin{itemize}
	\item   prediction horizon $\Np$ and control horizon $\Nc$;
	\item  positive-semidefinite weight matrices $Q_{y}$, $Q_{u}$,  $Q_{\Delta u}$; 
	\item positive constant $Q_{\varepsilon}$      used to soften the  constraints and thus to guarantee feasibility of the optimization problem~\eqref{eq:MPC};
	\item tolerances used in the QP algorithm solving~\eqref{eq:MPC}--\eqref{eq:MPC_constraints}.
\end{itemize}


\noindent In order to compact the notation, all the MPC knobs are collected in a  parameter vector $\parvec \in \setparvec \subseteq \mathbb{R}^{\sizeparvec}$, where $\setparvec$ is a given bounded set in which the optimal tuning for $\parvec$ is sought.

In this paper, we present an experiment-driven approach  to tune the MPC knobs $\parvec$ in order to optimize the overall closed-loop performance based on pairwise preferences expressed by a controller calibrator.  More precisely, we  assume that we do not have a  quantitative scoring function  used by the calibrator to measure  the closed-loop performance index and thus to be usable for auto-tuning $\parvec$ by global optimization. The only assumption we make is that for a  given pair of different MPC parameters   $\parvec_1,\parvec_2\in \setparvec$, only a preference  expressed by the calibrator is available in terms of ``performance achieved with parameters $\parvec_1$ was better (or worse, or similar) than the one achieved with $\parvec_2$''.  
 
The rationale of our problem formulation is that, in multi-criteria decision making   such as in control system design, it is difficult (and sometimes impossible)   to quantify an overall scoring function, but anyway it is easier   for a calibrator to assess a  \emph{preference} between the outcomes of two experiments.

 In order to compute the optimal MPC parameters $\parvec$, the active preference learning  algorithm  GLISp~\cite{BePi20} is used, which proposes the experiments to be performed for  tuning   $\parvec$, through pairwise comparisons, that we review in the next section.

\section{Preference-based tuning}\label{sec:pre_tunning}

We adopt the active preference learning algorithm GLISp proposed in~\cite{BePi20}
to iteratively suggest a sequence of MPC parameters $\parvec_1,\ldots,\parvec_N$ to be tested and compared such that $\parvec_N$ approaches the ``optimal'' 
combination of parameters as $N$ grows.

Formally, given two possible MPC knobs $\parvec_1$ and $\parvec_2$, we define  
the \emph{preference function}  $\pref:\rr^{\sizeparvec}\times\rr^{\sizeparvec}\to\{-1,0,1\}$  as
\begin{equation}
\pref(\parvec_1,\parvec_2)=\left\{\ba{ll}
-1 &\mbox{if $\parvec_1$ ``better'' than $\parvec_2$}\\
0 &\mbox{if $\parvec_1$ ``as good as'' $\parvec_2$}\\
1 &\mbox{if $\parvec_2$ ``better'' than $\parvec_1$}.
\ea\right.
\label{eq:pref_fun}
\end{equation}
The modeling assumption behind GLISp is that the calibrator assigns preferences based on an underlying closed-loop performance index $J$ that he or she wants to minimize:
\begin{equation}
\pref(\parvec_1,\parvec_2)=\left\{\ba{ll}
-1 &\mbox{if $  J(\parvec_1)<  J(\parvec_2)$}\\
0 &\mbox{if $  J(\parvec_1)=  J(\parvec_2)$}\\
1 &\mbox{if $  J(\parvec_1)>  J(\parvec_2)$},
\ea\right.
\label{eq:pref_fun-f}
\end{equation}
where $J$ is totally unknown (to the algorithm and to the calibrator as well).
Only the evaluations $\pi(\parvec_1,\parvec_2)$ of the preference function are accessible.  

The goal of the preference-based MPC calibration approach  discussed in this paper is  to find the optimal MPC parameters $\parvec^{\star} \in \setparvec$ such that $\parvec^{\star}$
is ``better'' (or ``no worse'') than any other parameter $\parvec$ according to the preference function $\pref$, namely
\begin{equation}
\parvec^\star\ \mbox{such that}\ \pref(\parvec^\star,\parvec)\leq 0,\ \forall \parvec \in \setparvec. 
\label{eq:glob-opt-pref}
\end{equation} 
If such a vector $\parvec^\star$ is found, as~\eqref{eq:glob-opt-pref} implies that $J(\parvec^\star)\leq J(\parvec)$,
$\forall \parvec \in \setparvec$, clearly we have also found a global minimizer  $\parvec^\star$ of $J$ on $\setparvec$.

\subsection{Training a surrogate function from preferences}

Assume that we have generated $N\geq 2$ samples $\{\parvec_1\ \ldots\ \parvec_N\}$ of the decision vector, 
with $\parvec_i,\parvec_j\in\rr^{\sizeparvec}$ such that $\parvec_i\neq \parvec_j$, $\forall i\neq j$,
$i,j=1,\ldots,N$. For each of these parameters, a closed-loop experiment is performed and the calibrator has provided  a \emph{preference vector} $B=[b_1\ \ldots\ b_{M}]^{T}\in\{-1,0,1\}^{M}$ with 
\begin{equation}
b_h=\pref(\parvec_{i(h)},\parvec_{j(h)}),
\label{eq:pref-vector}
\end{equation}
where  $M$ is the number of expressed preferences,  $1 \leq M \leq \binom{N}{2}$, 
$h\in\{1,\ldots,M\}$, $i(h), j(h)\in \{1,\ldots,N\}$, $i(h) \neq j(h)$. Note that the element $b_h$ of vector $B$  represents the preference expressed by the calibrator between the  closed-loop performance achieved with parameters  $\parvec_{i(h)}$ and $\parvec_{j(h)}$.

The observed preferences are then used to learn a surrogate function $\hat J:\rr^{\sizeparvec} \to\rr$ of the underlying performance index $J$.  The surrogate $\hat J$ is constructed by imposing the constraints
\begin{equation}
\hat\pref(\parvec_{i(h)},\parvec_{j(h)})=\pref(\parvec_{i(h)},\parvec_{j(h)}),\ \forall h=1,\ldots,M,
\label{eq:surrogate-condition}
\end{equation}
on $\hat J$, where $\hat \pref$ is defined from $\hat J$ as in~\eqref{eq:pref_fun-f}. 

The  function  $\hat J$ is parametrized as the following linear combination of
Radial Basis Functions (RBFs)~\cite{Gut01,MGTM07}: 
\begin{equation}
\hat J(\parvec)=\sum_{k=1}^N\beta_k\phi(\epsilon d(\parvec,\parvec_i)),
\label{eq:rbf}
\end{equation}
where  $d:\rr^{\sizeparvec}\times\rr^{\sizeparvec}\to\rr$ is the squared
Euclidean distance
\begin{equation}
d(\parvec_1,\parvec_2)=\|\parvec_1-\parvec_2\|_2^2,
\label{eq:distance}
\end{equation}
$\epsilon>0$ is a scalar parameter, $\phi:\rr\to\rr$ is an RBF,
and $\beta=[\beta_1\ \ldots\ \beta_N]^{T}$ are the unknown coefficients to be computed based on the preferences imposed in~\eqref{eq:surrogate-condition}. 
Examples of RBFs are $\phi(\epsilon d)=\frac{1}{1+(\epsilon d)^2}$
(\emph{inverse quadratic}), $\phi(\epsilon d)=e^{-(\epsilon d)^2}$ (\emph{Gaussian}),
$\phi(\epsilon d)=(\epsilon d)^2\log(\epsilon d)$ (\emph{thin plate spline}), see
more examples in~\cite{Gut01,Bem20}.

According to~\eqref{eq:surrogate-condition} and the preference relation~\eqref{eq:pref_fun-f}, the following constraints are imposed
on $\hat J$: 
\begin{equation}
\ba{ll}
\hat J(\parvec_{i(h)})\leq\hat J(\parvec_{j(h)})-\sigma+\varepsilon_h& \mbox{if}\ \pref(\parvec_{i(h)},\parvec_{j(h)})=-1\\
\hat J(\parvec_{i(h)})\geq\hat J(\parvec_{j(h)})+\sigma-\varepsilon_h& \mbox{if}\ \pref(\parvec_{i(h)},\parvec_{j(h)})=1\\
|\hat J(\parvec_{i(h)})-\hat J(\parvec_{j(h)})|\leq \sigma+\varepsilon_h& \mbox{if}\ \pref(\parvec_{i(h)},\parvec_{j(h)})=0
\ea
\label{eq:RBF-pref}
\end{equation}
for all $h=1,\ldots,M$, where $\sigma>0$ is a given tolerance
and $\varepsilon_h$ are positive slack variables. 

Accordingly, and similarly to Support Vector Machines (SVMs)~\cite{SS04}, the coefficient vector $\beta$ describing the surrogate $\hat J$  is obtained by solving the convex  QP problem
\begin{equation}
\ba{rll}
\min_{\beta,\varepsilon} &\displaystyle{\sum_{h=1}^{M} c_h\varepsilon_h+\frac{\lambda}{2}\sum_{k=1}^N\beta_{k}^2 }\\
\st 
& \displaystyle{\sum_{k=1}^N(\phi(\epsilon d(\parvec_{i(h)},\parvec_{k})-\phi(\epsilon d(\parvec_{j(h)},\parvec_{k}))\beta_k}
\\&\quad\quad\leq -\sigma+\varepsilon_h,\hfill \forall h:\ b_h=-1\\
&\displaystyle{\sum_{k=1}^N(\phi(\epsilon d(\parvec_{i(h)},\parvec_{k})-\phi(\epsilon d(\parvec_{j(h)},\parvec_{k}))\beta_k}\\
&\quad\quad\geq \sigma-\varepsilon_h, \hfill \forall h:\ b_h=1\\
&\displaystyle{\left| \sum_{k=1}^N(\phi(\epsilon d(\parvec_{i(h)},\parvec_{k})-\phi(\epsilon d(\parvec_{j(h)},\parvec_{k}))\beta_k\right|}\\&
\quad\quad\leq\sigma+\varepsilon_h,\hfill \forall h:\ b_h=0\\  & h=1,\ldots,M
\ea
\label{eq:QP-pref}
\end{equation}
where $c_h$ are positive weights, for example $c_h=1$, $\forall h=1,\ldots,M$.
The slack variables $\varepsilon_h$ in~\eqref{eq:QP-pref}
are used to relax the constraints imposed by  the preference vector $B$
in~\eqref{eq:surrogate-condition}.
Constraint infeasibility might be due to an inappropriate selection
of the RBF (namely, poor flexibility in the parametric description of the surrogate $\hat J$) and/or to outliers in the acquired preferences, for instance
due to inconsistent assessments done by the calibrator. The scalar $\lambda$ in the cost function~\eqref{eq:QP-pref} is a regularization parameter. When $\lambda>0$, problem~\eqref{eq:QP-pref} is a QP problem that admits a unique solution because 
$c_h>0$ for all $h=1,\ldots,M$. If $\lambda=0$, problem~\eqref{eq:QP-pref} becomes a Linear Program (LP), whose solution may not be unique.

We remark that computing the surrogate function $\hat J$  requires one to choose the hyper-parameter $\epsilon$ defining the shape of the RBFs $\phi$ in~\eqref{eq:rbf}. This parameter can be   chosen through $K$-fold cross-validation~\cite{stone1974cross}, by testing  the capabilities of  $\hat J$ in reconstructing the preferences in slices  of the dataset not used to estimate $\hat J$.


\subsection{Acquisition function}
\label{sec:acquisition}
Once a  surrogate  $\hat J$ is estimated, this function can be in principle minimized in order to find the optimal MPC parameter vector $\parvec$. More specifically, the following steps can be followed:  
($i$) generate a new sample by pure minimization of the estimated surrogate function $\hat J$ 
defined in~\eqref{eq:rbf}, \emph{i.e.}, 
\[
\parvec_{N+1}=\arg\min\hat J(\parvec)\ \st\   \parvec\in\setparvec; 
\] ($ii$) ask the calibrator to evaluate the preference
$\pref(\parvec_{N+1},\parvec^\star_N)$, where $\parvec^\star_N\in\rr^{\sizeparvec}$ is the best vector of MPC parameters found so far, corresponding
to the smallest index $i^\star$ such that
\begin{equation}
\pref(\parvec_{i^\star},\parvec_i)\leq 0,\ \forall i=1,\ldots,N;
\label{eq:pref-x-star}
\end{equation}
 ($iii$) update the estimate of $\hat J$ through~\eqref{eq:QP-pref}; and ($iv$) iterate over $N$.

Such a procedure, which only \emph{exploits} the current available observations in finding the optimal MPC parameter vector $\parvec$, may easily miss the global minimum of~\eqref{eq:glob-opt-pref}. Therefore, looking only at the surrogate function $\hat J$ is not enough to search for a new sample $\parvec_{N+1}$. A term promoting the \emph{exploration} of the parameter space should thus be considered. 

In the GLISp algorithm, an acquisition function is employed to balance exploitation \emph{vs.}  exploration when generating the new sample $\parvec_{N+1}$. 
As proposed in~\cite{Bem20}, the exploration function is constructed by   using the inverse distance weighting (IDW) function $z:\rr^{\sizeparvec}\to\rr$ defined by
\begin{equation}
z(\parvec)=\left\{\ba{ll}
0 & \mbox{if}\ \parvec\in\{\parvec_1,\ldots,\parvec_N\}\\
\tan^{-1}\left(\frac{1}{\sum_{i=1}^Nw_i(\parvec)}\right)&
\mbox{otherwise}\ea\right.
\label{eq:IDW-distance}
\end{equation}
where $w_i(\parvec)=\frac{1}{d^2(\parvec,\parvec_i)}$. Clearly $z(\parvec)=0$ for all parameters already tested, and $z(\parvec)>0$ in $\rr^{\sizeparvec}\setminus \{\parvec_1,\ldots,\parvec_N\}$.  The arc tangent function in~\eqref{eq:IDW-distance} avoids that $z(\parvec)$ gets excessively large far away from all sampled points.




Then, given an exploration parameter $\delta\geq 0$, the \emph{acquisition function} $a:\rr^{\sizeparvec}\to\rr$
is constructed as
\begin{equation}
a(\parvec)=\frac{\hat J(\parvec)}{\Delta\hat J}-\delta z(\parvec),
\label{eq:acquisition}
\end{equation}
where 
\[
\Delta \hat J=\max_i\{\hat J(\parvec_i)\}-\min_i\{\hat J(\parvec_i)\}
\]
is   the range of the surrogate function on the samples in $\{\parvec_1,\ldots,\parvec_N\}$ and is used in~\eqref{eq:acquisition} as a normalization factor to simplify the choice of the exploration
parameter $\delta$.  
Clearly $\Delta \hat J \geq \sigma$ if at least one comparison $b_h=\pref(\parvec_{i(h)},\parvec_{j(h)}) \neq 0$.  

As  discussed below, given a set $\{\parvec_1,\ldots,\parvec_N\}$ of samples and a vector $B$
of preferences defined by~\eqref{eq:pref-vector}, the next MPC parameter  $\parvec_{N+1}$ to  test  
is computed as the solution of  the (non-convex) optimization problem
\begin{equation}
\parvec_{N+1}=\arg\min_{\parvec\in \setparvec
} a(\parvec).
\label{eq:xNp1}
\end{equation}



The preference-based global optimization algorithm based on
the RBF interpolants~\eqref{eq:rbf}  and the acquisition function~\eqref{eq:acquisition} is detailed in~\cite{BePi20}
and will be referred to as the GLISp algorithm since now on.
We will use GLISp to compute the optimal MPC parameter vector $\parvec^\star$ by preferences.
We use the Particle Swarm Optimization (PSO) algorithm of~\cite{VV09} to solve problem~\eqref{eq:xNp1}. Note that the construction of the acquisition function $a$ is rather heuristic, therefore finding highly accurate  solutions of~\eqref{eq:xNp1}
is not required. 

In the acquisition function~\eqref{eq:acquisition}, the exploration parameter $\delta$ promotes sampling the space in $\setparvec$ in areas that have not been explored yet. While, as observed earlier, $\delta=0$ can make the GLISp algorithm
rely only on the surrogate function $\hat J$ and miss the global optimum,
setting $\delta\gg1$ makes the GLISp algorithm exploring the entire feasible region regardless of the results of the comparisons. For a sensitivity analysis example
of GLISp with respect to $\delta$, the reader is referred to~\cite[Section 7.5]{BePi20}.
Regarding the other main hyper-parameter $\epsilon$ of GLISp defining the RBF in~\eqref{eq:rbf}, during the active learning phase  $K$-fold cross-validation is executed repeatedly to automatically choose and possibly adapt $\epsilon$.

We finally remark that, in executing the GLISp, some values of the MPC knobs may lead to an unstable closed-loop behaviour. In this case, the experiment should be interrupted (e.g., for safety reasons) and the calibrator may still express a preference. It is also possible that in the initialization phase of GLISp a comparison should be performed over two experiments both exhibiting closed-loop instability. In this case,  one experiment can be preferred to another one,  for instance, based on  the time it could run before being interrupted.
 


\section{Cases studies}\label{sec:case_study}
We apply the proposed preference-based MPC calibration method on two case studies: control of a \emph{Continuous Stirring Tank Reactor}  (CSTR) and autonomous driving with obstacle avoidance. The role of the calibrator is played by  the first author of the paper, which compares the observed closed-loop performance manually  according to a qualitative  scoring function constructed  in her mind based on engineering insights.

In both case studies we set the maximum number of function evaluations $N_{\rm max}=50$, the exploration parameter $\delta=0.3$
, and the tolerance parameter $\sigma=10^{-6}$. 
The GLISp algorithm is initialized with
$N_{\rm init}=10$ random samples generated by Latin hypercube sampling~\cite{MBC79} using the \emph{lhsdesign} function of  the  \emph{Statistics  and  Machine  Learning Toolbox} of MATLAB ~\cite{matlabStaML}. The  hyper-parameter $\epsilon$ defining the   RBF functions~\eqref{eq:rbf} is initialized as $1$ and updated at iterations $10,20,30,40$ via $K$-fold cross-validation with $K=3$. The tests are run on  an x64 with Intel i7-8550U 1.8 GHz CPU and 8GB of RAM. The \emph{Model Predictive Control Toolbox} of MATLAB~\cite{matlabMPC} is used for MPC design and simulation.


\subsection{CSTR optimal steady-state switching policy}
\label{sec:CSTR}
\subsection*{\underline{System description}}
 The first case study considers the control of a CSTR, extensively described  in~\cite{bequette1998process}. The CSTR system  consists of  a jacketed non-adiabatic tank where an exothermic reaction  occurs.  The tank is assumed to be perfectly mixed with constant inlet and outlet rate. 
The chemical laws describing the CSTR reaction can be derived based on the energy and material balance, and are given by:
\begin{align}
 &\frac{dT(t)}{dt} \!  \ \ =  \frac{F}{V}(T_f(t) \! -\!  T(t))\!  -  \! \frac{H}{C_p \rho} r(t) \! - \!  \frac{US}{C_p \rho  V}(T(t)\! -\! T_c(t))\nonumber\\
  &\frac{dC_A(t)}{dt} \!  =   \frac{F}{V}(C_{Af}(t)\!  - \! C_A(t)) \! -\!  r(t), \label{cstr_system_eqn}
\end{align}
where $V$ [m$^3$] is the volume of the reactor, and $F$ [m$^3$/hr] is the rate of reactant $A$ feeding the tank, which is equal to the rate of the product stream that exists from the reactor. 
Both $V$ and $F$ are assumed to be constant. $C_A(t)$ and $C_{A_f}(t)$ [kgmol/m$^3$] represent  the concentration at time $t$ of reactant $A$ in the tank and in the  inlet feed stream, respectively.  $T(t)$, $T_f(t)$ and $T_c(t)$ [K] are the temperature of the reactor, of the inlet feed stream and of the coolant stream, respectively. Constant heat of reaction $H$ [kcal/kgmol], fluid heat capacity $C_p$ [kcal/(kg K)] and density $\rho$ [kg/m$^3$] are assumed. $U$ [kcal/(m K hr)] and $S$ [m$^2$] are the overall heat transfer coefficient and heat exchange area, respectively. $r(t)$ [kgmol/(m$^3$hr)] is the reaction rate per unit volume, which can be calculated through Arrhenius rate law: 
\begin{equation}\label{eq:ratelaw}
r(t) = k_0 \exp\left(\frac{-E}{R T(t)}\right)\ C_A(t)
\end{equation}
where $k_0$ [hr$^{-1}$] is the pre-exponential factor; $E$ \mbox{[kcal/kgmol]} is the activation energy for the reaction; and $R$ \mbox{[Kcal/(kgmol K)]} is the gas constant. The CSTR process parameters are taken from \cite{bequette1998process} and listed in Table~\ref{tab:CSTR_para}.

\begin{table}[!bt]
	\caption{CSTR parameters
	 \cite{bequette1998process}}\label{tab:CSTR_para}
	\centering
	\vspace{-0.2cm}
	\begin{tabular}{lllll}
		\hline
		Parameter & Value      & unit \\
		\hline
		$F/V$     & 1.0          & $\text{hr}^{-1}$  \\
		$k_0$    &  3.49$e7$    & $\text{hr}^{-1}$      \\
		$H$      & 5960       & $\text{kcal}/\text{kgmol}$  \\
		$E$      & 11843      & $\text{kcal}/\text{kgmol}$  \\
		$\rho C_p$      & 500      & $ \text{kcal}/(\text{m}^3 \text{K})$  \\
		$US/V$      & 150      & $ \text{kcal}/(\text{m}^3 \text{K hr})$  \\
		$T_{c0}$     & 298          & $\text{K}$  \\
		$R$     & 1.987          & $\text{Kcal/(kgmol K)}$  \\
		\hline
	\end{tabular}
	\vspace{-0.4cm}
\end{table}

 \subsection*{\underline{Control objectives}}
Initially, the plant is operating at  steady state with a reactant concentration  $C_A = 8.56 $ kgmol/m$^3$ and  low conversion. The objective is to design an MPC
controller to achieve a steady-state concentration $C_A = 2$ kgmol/m$^3$, thus increasing the conversion rate. The feed stream concentration $C_{A_f}(t)$ and temperature $T_f(t)$ are treated as measured disturbances, and the coolant temperature  $T_c(t)$ is  the control input. In this test, the condition of the feed stream is kept constant, with $T_f$ = 298.15 K and $C_{A_f}$ = 10 kgmol/m$^3$.

At  each time step, the nonlinear model  of the CSTR system in~\eqref{cstr_system_eqn} is linearized around its operating point, and a linear MPC problem is solved. Although the output variables of interest are both the reactor temperature  $T(t)$ and the concentration $C_A(t)$ in the product stream,  only the latter is tracked and  compensated by MPC.



Different competitive objectives should be taken into account in the MPC calibration, such as fast steady-state transition, reasonable final target achievement, and low energy consumption.
The following  guidelines are used to assist the calibration. First, the steady-state switching is expected to be completed within two days.  Second, the final achieved  steady state should be within $\pm3$\% of the desired value $C_A = 2$ kgmol/m$^3$. Third, in order to take into account  energy consumption due to the cooling process,  the temperature of the coolant stream  $T_c(t)$ is restricted to be in the range of $[284 \ \ 310]$ K, with a maximum change at each time step ($T_{c_{max}}$) set to $10$ K. The first two requirements just reflect the desired performance, and thus  are not treated as  \emph{hard} constraints during calibration.  

 The following three design parameters are tuned:  the sampling time $T_s$ used to close the loop and to discretize the continuous-time model~\eqref{cstr_system_eqn}; the prediction horizon $N_p$ and the weight $Q_{\Delta u}$ penalizing the input change. These tuning parameters $T_s$, $N_p$ and $\log(Q_{\Delta u})$ are restricted to the ranges $[0.25\  1.5]$ hr, $[4 \ 40]$ and $[-5 \  3]$, respectively. The control horizon  $N_u$ is set equal to  $N_p/3$ and rounded to the closet integer. Other MPC  knobs are fixed, with $Q_y$ and $Q_u$ equal to $\left[\begin{smallmatrix}0 & 0 \\0 & 1 \end{smallmatrix}\right]$ and 0, respectively. 

 \subsection*{\underline{Calibration process}}

At each iteration of GLISp, a set of new MPC parameters is suggested and the closed-loop experiment is performed in simulation. The experiment is  interrupted either when the  steady-state condition is reached or after $48$ hrs (namely, two simulated days). The experiment is also interrupted if an un-safe behavior is observed. The calibrator is then asked to compare the performance of the experiment with the best   performance achieved  until that time
and  choose the one she prefers based on the control objectives and the aforementioned guidelines.

\begin{figure}[t!]
  \centering
  \includegraphics[width=\linewidth]{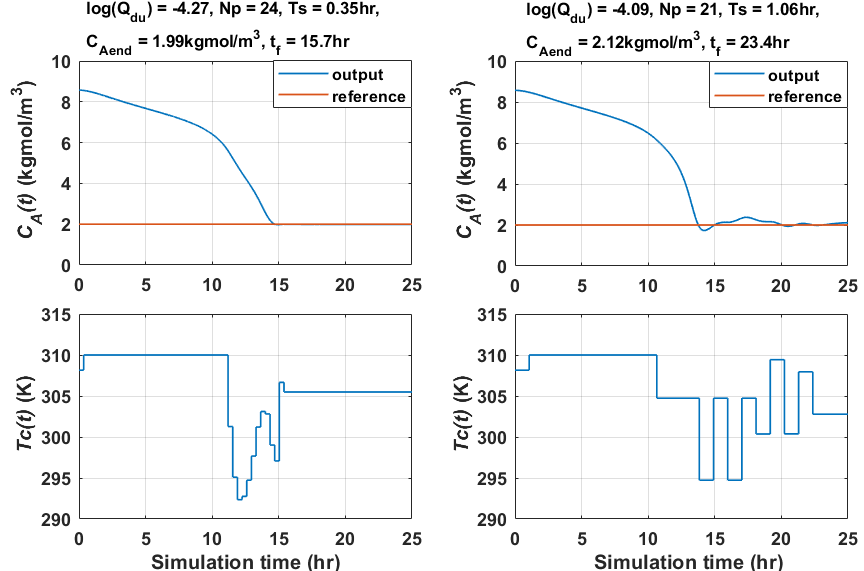}
  \vspace{-0.6cm}
  \caption{CSTR query window for one iteration of GLISp. Left experiment is preferred to the right one because of a faster transient time.}
  \label{fig:CSTR_query}
\end{figure}

As an example, Fig. \ref{fig:CSTR_query} shows one iteration of the semi-automated calibration process. At the top of the figure, the MPC design parameters, the achieved values of $C_{A_{end}}$ and duration of the switching process $t_f$ are displayed. Both the MPC designs achieved the desired $C_{A_{end}}$ within $48$ hrs. However, the transient and settling  times  of the  left-hand-side experiment are shorter than the right-hand-side one. Moreover,  the variations of the input signal $T_{c}$ are smaller for the experiment on the left. Therefore, in this iteration, the left-hand-side MPC design is preferred.

 \subsection*{\underline{Results}}

The GLISp algorithm terminates after $N_{\max}=50$ iterations (namely, closed-loop experiments) and hence $49$ comparisons.   The best MPC design parameters $T_s$, $N_p$ and $\log(Q_{\Delta u})$ are found to be $0.31$ hr, $26$ and $-1.79$, respectively.  Fig. \ref{fig:CSTR_final_performance} (left panels) shows the corresponding trajectory of reactant concentration  $C_A(t)$ and of the  manipulated variable $T_c(t)$.

\subsection*{\underline{Fully-automated calibration}}

For the sake of comparison, a fully-automated  calibration is performed. This requires to define a multi-objective quantitative scoring function describing the expected closed-loop performance. This step can be very hard and time-consuming, requiring proper scaling and balancing of the competitive control objectives.  

The following closed-loop performance index function to be minimized is constructed after several trial-and-error iterations based on the preferences of the calibrator:
\begin{equation}\label{eq:perf_index}
\begin{split}
  \frac{t_f}{t_{f_{max}}} \! + \!  \frac{\sum_{k=1}^{N_T} (T_{c_k} - T_{c_{k-1}})^{2}}{T_{c_{max}} \ N_T} \! + \!  \frac{\left|   C_{A_{end}}-C_{A_{ref}}  \right|}{AR\% \ C_{A_{ref}}} 
\end{split}
\end{equation}
where  $t_{f_{max}}=48$ hrs is the  maximum duration of the switching process; $T_{c_{max}}=10$ K is the maximum allowed temperature change of the coolant fluid between each time step; $T_{c_k}$ and $T_{c_{k-1}}$ are  the temperature of the coolant fluid at time step $k$ and $k-1$, respectively; $N_T$ is the total number of time steps in the steady-state transition; 
$AR\% = 3\%$ is the acceptable range of the final steady state concentration; $C_{A_{end}}$ and $C_{A_{ref}}$ [kgmol/m$^3$] are the achieved and desired final steady-state concentration of reactant $A$.  

The GLIS algorithm (without preference)  of~\cite{Bem20} is used for experiment-driven optimization of the cost function in~\eqref{eq:perf_index} with respect to the MPC design parameters. For a fair comparison with the proposed preference-based algorithm, $50$ function evaluations are performed. This value actually does not account for the number of trials (namely, experiments) needed to construct  the scoring function  in~\eqref{eq:perf_index}.

\begin{figure}[t!]
  \centering
  \includegraphics[width=\linewidth]{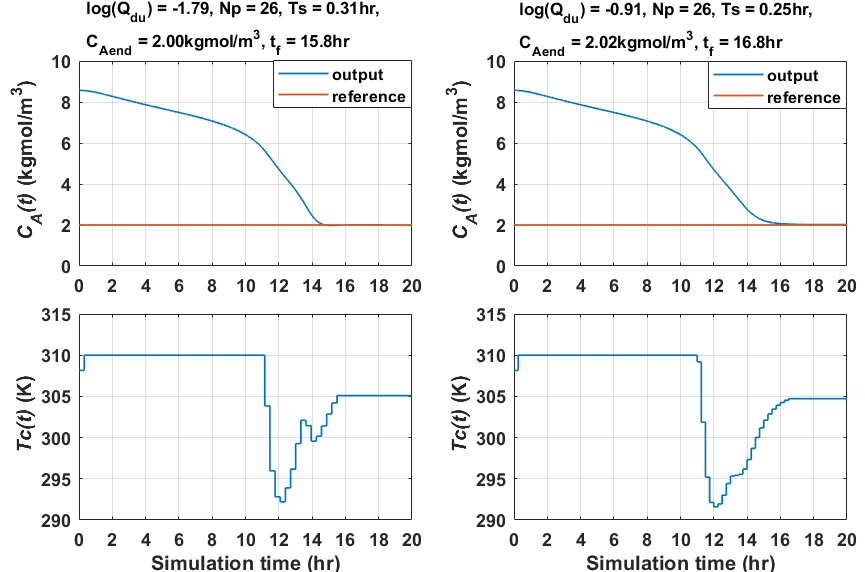}
  \vspace{-0.6cm}
  \caption{CSTR closed-loop performance obtained by calibrating MPC parameters through the proposed semi-automated preference-based approach (left panels)   and through a fully-automated approach minimizing the scoring function~\eqref{eq:perf_index} (right panels).}
  \label{fig:CSTR_final_performance}
\end{figure}

The obtained optimal MPC parameters $T_s$, $N_p$ and $\log(Q_{\Delta u})$ are $0.25$ hr, $26$ and $-0.91$. Fig. \ref{fig:CSTR_final_performance} (right panels) shows the corresponding closed-loop trajectory of $C_A(t)$ and of the  manipulated variable $T_c(t)$. It can be observed that the preference- and the non-preference- based approaches achieve similar closed-loop performance. 

Overall, this case study demonstrates the capability of the semi-automated approach for solving calibration tasks with multiple competitive objectives. For such calibration tasks, when using a fully-automated approach, significant efforts need to be devoted to construct a proper performance scoring function as in~\eqref{eq:perf_index}. Therefore, the semi-automated approach can greatly reduce calibration time and efforts by eliminating this step. 




 
\subsection{Autonomous driving vehicle} 

 \subsection*{\underline{System description}}
As a second case study, we consider the problem of  lane-keeping (LK) and obstacle-avoidance (OA) in  autonomous driving. A simplified two degree-of-freedom bicycle model   is used to describe the vehicle dynamics,  with the front wheel as the reference point. The model involves the following three state variables:  longitudinal  $x_f$ and lateral  $y_f$ [m] position of the front wheel, and  yaw angle $\theta$ [rad].  The control inputs are   the reference for the vehicle velocity $v$ [m/s] and the steering angle $\delta_s$ [rad]. The system behaviour is thus described by the continuous-time kinematic model 
\begin{equation}\label{eq:nonlin_bicyclemodel}
\begin{split}
 \dot{x}_f = & v \cos(\theta+\delta_s)\\
 \dot{y}_f = & v \sin(\theta+\delta_s)\\
 \dot{\theta} = & \frac{v\sin(\delta_s)}{L}
\end{split}
\end{equation}
where $L$ [m] is the vehicle length. This model is  linearized around its operating point at each time step and used to design a linear MPC. 

 \subsection*{\underline{Control objectives}}
In tuning the MPC parameters,  the  objective is to keep the vehicle at the same horizontal lane with constant speed and to overtake other moving vehicles in an optimal way if they are within safety distance. However, similar to the CSTR case, it is difficult to define a proper quantitative  scoring function for this multi-objective calibration task (for example, due to the ambiguity of transferring ``optimal obstacle avoidance'' into a mathematical formula).  Therefore, to achieve good MPC performance using a fully-automated approach not based on preferences, either the calibrator needs to find a proper scoring function via trial-and-error (like in the previous CSTR case)  or an advanced path planner model is required to provide a well-defined reference path to compare with. Both procedures can be time-consuming and computationally heavy. On the other hand, when using the semi-automated approach discussed in this paper, none of the pre-mentioned steps is required.  

The test scenario is the following. The vehicle starts at position ($x_f$, $y_f$) = (0, 0) m with $\theta$ equal to 0 rad and there is another vehicle (obstacle) at position ($30$, $0$) m  which is moving at constant speed ($40$ km/hr), with a constant yaw angle of $0$ rad.  Both vehicles are assumed to have a simplified rectangular shape with length equal to $4.5$ m and width of $1.8$ m. At nominal LK condition, the vehicle moves at $50$ km/hr with $y_f$ = 0 \mbox{ m}. When the obstacle is within safety distance ($10$ m in this case), the vehicle needs to overtake it and the minimum lateral distance between two vehicles is set to $3$ m. During LK and OA periods, the vehicle is allowed to vary its velocity in the range of $[40 \  70]$ and $[50\  70]$  km/hr, respectively, with its reference velocity set to $50$ and $60$ km/hr, respectively.   The maximum and minimum of $\theta$ in both LK and OA periods are $\pm \frac{\pi}{4}$ rad with the maximum and  minimum rate of change between each time step set to $\pm 0.0873$ rad/s.

Five MPC design parameters are  tuned. The sampling time $T_s$ is allowed to vary in the range $[0.085 \  0.5]$ s. The prediction horizon $N_p$ is restricted to $[10 \  30]$ and the control horizon $N_u$ is taken as a fraction $\epsilon_c$ of $N_p$ rounded to the closest integer. Here, $\epsilon_c$ can take values in the range $[0.1 \  1]$. The weight matrix of manipulated variables ($Q_{\Delta u}$)  is set to be diagonal with diagonal entries $q_{u11}$ and $q_{u22}$ such that  $\log(q_{u11})$ and $\log(q_{u22})$ vary in the range $[-5 \  3]$. The other MPC design parameters are fixed, with $Q_y$ and $Q_u$  set to $\left[\begin{smallmatrix}1 & 0 & 0 \\0 & 1 & 0\\0 & 0 & 0 \end{smallmatrix}\right]$ and $\left[\begin{smallmatrix}0 & 0 \\0 & 0 \end{smallmatrix}\right]$, respectively. 

 \subsection*{\underline{Calibration process}}

The calibrator  selects the preferred controller based on the following observations: ($i$) during both LK and OA periods, the worst-case computational time ($t_{comp})$ required for solving the QP problem~\eqref{eq:MPC} at each time step  needs to be smaller than $T_s$, so that the  MPC can be implemented in real time; ($ii$) during the LK phase, the vehicle should move at constant speed with $y_f$ and $\theta$ close to $0$ m and $0^\circ$; ($iii$) during the OA phase, the vehicle should keep reasonable safety distance away from the obstacle and should guarantee passengers' comfort  (\emph{i.e.},  aggressive lateral movements during overtaking should be avoided); ($iv$) the velocity in both LK and OA period should be close to the reference value with its variations kept to minimum; ($v$) variations of steering angles should not be aggressive; ($vi$) when there is a conflict combination among aforementioned criteria, criterion ($i$) has the highest priority and if the conflict is among criteria ($ii$)--($v$), preference is given to the one leading to safer driving practice based on calibrator's experience.

\begin{figure}[t!]
  \centering
  \includegraphics[width=\linewidth]{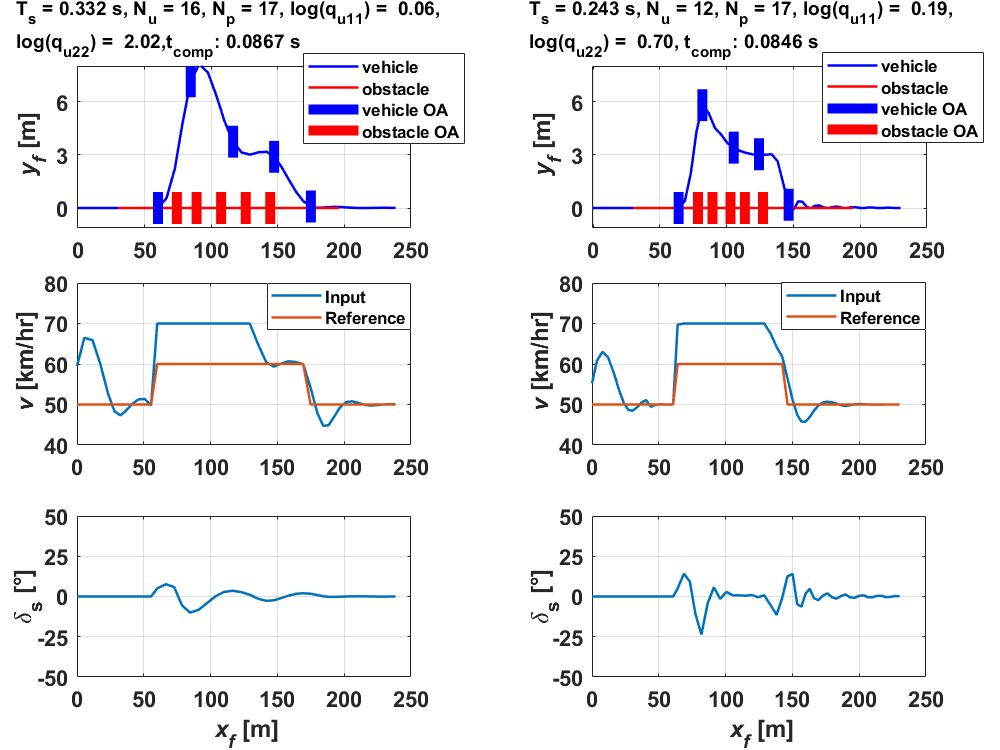}
  \vspace{-0.7cm}
  \caption{Vehicle control query window. The top subplots show the location trajectories of the vehicle and the obstacle, in which the ``vehicle OA'' and ``obstacle OA'' bars show five relative positions of the vehicle and obstacle during the obstacle-avoidance phase. The middle subplots show the actual and reference velocity $v$ at different longitudinal positions. The steering angle $\delta_s$ over the longitudinal position  is depicted in the bottom subplots. For ease of assessment, the unit of $v$ and $\delta_s$ in the figure is converted to \mbox{km/hr} and degree ($^\circ$), respectively. 
  	The results on the right panels are preferred. }
  \label{fig:car_comparison}
\end{figure}

The closed-loop test is simulated for 15 seconds. Fig. \ref{fig:car_comparison} shows the query window for one iteration of the calibration process. The MPC design parameters and $t_{comp}$ are displayed at the top of the figure. In both of the cases illustrated in the figure, the QP problem~\eqref{eq:MPC} is solved within the chosen sampling time $T_s$. The performance of the right-hand-side experiment is preferred since in the left-hand-side experiment  large lateral movements are present. Indeed, these movements  can be much more dangerous (as the car may cross to the other lane) comparing to slightly  more aggressive $\delta_s$ variations.


 \subsection*{\underline{Results}}
The GLISp algorithm terminates after $N_{\max}=50$  function evaluations and $49$ comparisons. The best MPC design parameters $T_s$, $\epsilon_c$, $N_p$, $\log(q_{u11})$ and $\log(q_{u22})$ are found to be equal to  $0.085$ s, $0.310$, $16$, $0.261$ and $0.918$, respectively with a worst-case computational time $t_{comp}$ needed to solve the QP problem~\eqref{eq:MPC} equal to $0.0808$ s. 
 
\begin{figure}[t!]
  \centering
  \includegraphics[width=\linewidth]{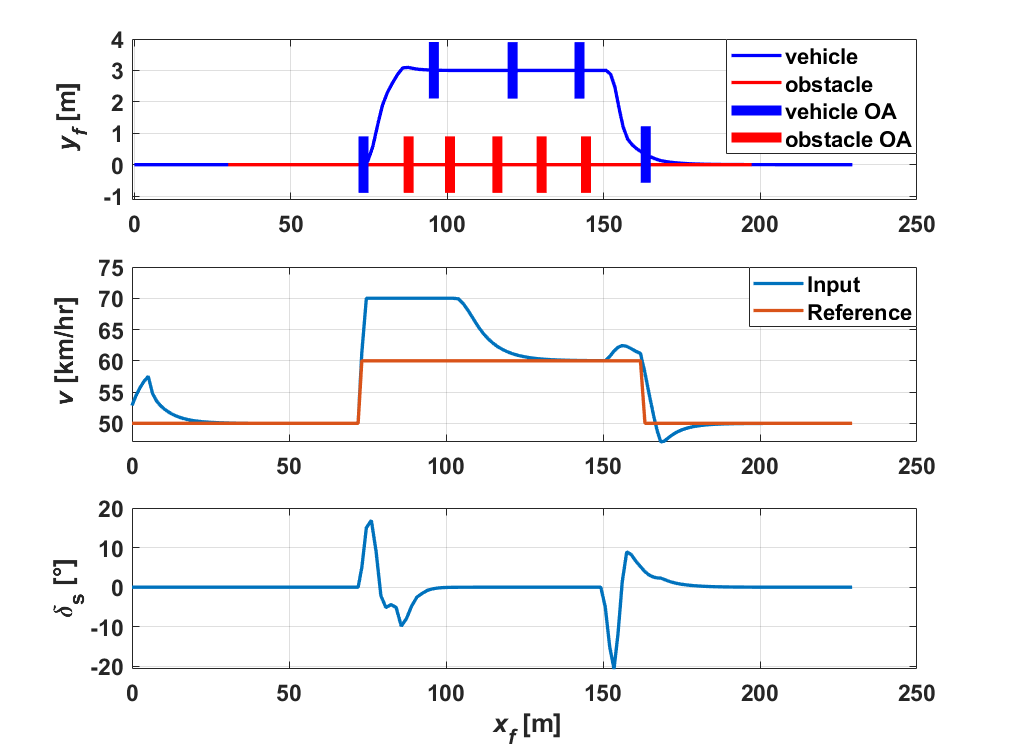}
  \vspace{-0.8cm}
  \caption{Vehicle control final performance obtained by the designed MPC controller. }
  \label{fig:car_finalPerf}
\end{figure}

The final closed-loop results   achieved by the designed MPC are depicted in Fig.~\ref{fig:car_finalPerf}, which demonstrates that solely based on calibrator's preferences, after only $50$ experiments, the proposed algorithm is able to tune the MPC parameters with  satisfactory performance. It is also important to remark that, for the same problem, the authors were not able to find, via several trial-and-error tests, a proper scoring function of  the closed-loop performance to be used for a fully-automated calibration.

\section{Conclusion} \label{sec:conclusion}
In this paper, a novel semi-automated MPC calibration approach was presented which turns out to be efficient in terms of number of experiments required for calibration. The key feature of the proposed methodology is that it allows calibration only based on pairwise preferences between the outcomes of the experiments, and thus it is very useful for calibration tasks with qualitative, subjective, or hard-to-quantify performance index functions, and with calibrators having limited MPC design knowledge.  The same preference-based  approach can  be also used for calibration of other type of controllers, such as PIDs.

Future research  activities are devoted to add new features to the calibration algorithm in order to further improve its performance. For instance,  when the outcome of the experiment can be easily labelled  as ``very bad'',    one can try to exploit this information instead of specifying a preference w.r.t. other experiments. A straightforward solution could be to train a classifier to separate the feasible space into good/bad regions, and thus reduce the search space.






\bibliographystyle{IEEEtran} 
\bibliography{Biblio_ECC}    



\end{document}